\documentclass[12pt,a4paper]{article}%
\usepackage{amsfonts}
\usepackage{amssymb}
\usepackage{epsfig}
\usepackage{graphics}
\usepackage{amsmath}
\usepackage[a4paper]{geometry}
\usepackage{fancyhdr}
\usepackage{amssymb}
\usepackage{color}
\usepackage{caption}
\usepackage{subcaption}
\usepackage{xfrac}
\usepackage{enumitem}
\usepackage[pdftex]{hyperref}
\usepackage{graphicx}%
\providecommand{\U}[1]{\protect\rule{.1in}{.1in}}
\hypersetup{                                                  a4paper,                                                  breaklinks                                            }

\renewcommand{\emph}[1]{\textsl{#1}}
\newlength{\bibitemsep}
\newlength{\bibparskip}
\let\oldthebibliography\thebibliography
\renewcommand\thebibliography[1]{%
  \oldthebibliography{#1}%
  \setlength{\parskip}{\bibitemsep}%
  \setlength{\itemsep}{\bibparskip}%
}
\begin{document}

\title{Computational Discovery with Newton Fractals, Bohemian Matrices, \& Mandelbrot Polynomials}
\author{%

\begin{tabular}
[l]{l}%
\textit{
Neil J.~Calkin, Eunice Y.~S.~Chan, and Robert M.~Corless$^{*}$}\\
calkin@clemson.edu, School of Mathematical and Statistical Sciences,\\ Clemson University, Clemson, SC 29634-0975\\
echan295@uwo.ca,
Centre for Medical Evidence, Decision Integrity and\\ Clinical Impact (MEDICI Centre), Department of Anesthesia and \\ Perioperative Medicine, Schulich School of Medicine and Dentistry,\\ Western University, London, Ontario, N6A 5B7 Canada\\
and $*$Corresponding Author: rcorless$@$uwo.ca,
School of Mathematical \\ and Statistical Sciences Western University, London, Ontario, N6A 5B7 Canada\\
\end{tabular}
}
\date{}
\maketitle

\begin{abstract}
The authors have been using a largely algebraic form of ``computational discovery'' in various undergraduate classes at their respective institutions for some decades now to teach pure mathematics, applied mathematics, and computational mathematics.  This paper describes what we mean by ``computational discovery,'' what good it does for the students, and some specific techniques that we used.
\end{abstract}


%


%
\thispagestyle{fancy}

\section{Setting the stage}
\begin{quote}
``The imparting of factual knowledge is for us a secondary consideration.  Above all we aim to promote in the reader a correct attitude, a certain discipline of thought, which would appear to be of even more essential importance in mathematics than in other scientific disciplines.'' \hfil P\'olya \& Szeg\H{o} vol I.~\cite[p.~VII]{polya1997problemsI}
\end{quote}

The preface quoted above from the classic book cited, which is nearly a hundred years old now, opens with an epigraph which we further paraphrase, as follows: ``What is good education? Giving students systematic opportunities to discover things for themselves.''
Indeed, \textsl{Computational Discovery}, also called ``Experimental Mathematics,'' is also very familiar to the research mathematician, not just mathematics educators: nearly everyone uses it (even if they say that they don't, or don't say that they do).  There can be no shame in it, if the likes of Gauss and Euler used the technique~\cite{borwein2004mathematics,Borwein2008}.  See also the excellent book~\cite{eilers2017introduction}. The most basic idea is, after all, very simple: one computes a few cases, tries to guess a pattern, and if successful, tries to prove it\footnote{In the hands of P\'olya, of course, this most basic idea was refined and extended almost beyond recognition. Similarly, Bill Bauldry advocates in~\cite{Bauldry:2020:ThreeAspects} a variation of this idea, which he terms the Action--Consequence--Reflection Principle; this is very similar to the advanced ideas that P\'olya advocates in~\cite{polya1962mathematical}. Still, it is astonishing how powerful this idea is, even at its simplest, given modern computational tools.}.  However, there is the all-too-common attitude that says ``Once one has a proof, who really cares how the theorem was discovered?''

Quite a few people care, actually.  First, it's much easier to get students to succeed in learning mathematics if we teach them this powerful method, which many mathematicians {actually} \emph{use}.
Second, the world of mathematics is larger than the human mind can normally conceive of: computational experiments can genuinely surprise even experienced mathematicians.  Third, mathematics is changing and has changed with the introduction of new powerful computational tools (for instance, the ordinary differential equation solvers in Julia are extremely impressive: see~\cite{Rackauckas2017}).  We have a duty to train our students in the tools used to explore that larger world.  We also have a duty to try to keep up with the changes in these tools: see for instance that aforementioned paper on the new solvers in Julia, but this has been known for a long while~\cite{abelson1976computation}.

Speaking of Julia, the phrase ``computational thinking'' has been given new life by an MIT course of that name, using Julia, offered by David Sanders, James Schloss, and Alan Edelman, and given away for free on~\href{https://youtu.be/vxjRWtWoD_w}{YouTube}.  Julia is sometimes referred to as a ``Python killer,'' because it is very nearly as easy to use as Python but also offers very nearly the same performance as C or Fortran.  Fortran, by the way, is very much still alive in the HPC world, for example in weather prediction.  Python is by no means dead yet either, and has several advantages for use in mathematics instruction.  Chief among these advantages is its very large user base amongst computational scientists: students are motivated to learn Python because it is directly used by computational neuroscientists, financial mathematicians (who may also use C), for machine learning, and many other buzzword-compliant research areas. The Open Educational Resource (OER) called ``Computational Discovery on Jupyter'' is intended to take advantage of these technologies~\cite{Chan:2022:CDJ}.

One user interface for Python, namely the Jupyter Notebook, seems to be dominating at present.  This allows access to Python, to Julia, to Matlab, and to SAGE Math, which offers an interface to many computer algebra systems.  A Jupyter Notebook allows combined use of text, formulas (through a flavour of ``Markdown'' which is kind of a blend of HTML and \LaTeX\ that is lighter weight than either), figures, data, and programs.  Those people who are used to Maple Documents (or Worksheets or Workbooks) will be able to see the similarity, and indeed Python is rather structurally similar to Maple as well.

An experienced experimental mathematician or instructor will know that the platform for teaching has some impact, both good and bad, on the quality of teaching.  One does not want the platform to \textsl{interfere} with the teaching.  One does not want to spend too much time in class explaining the platform, instead of using that valuable time explaining the concepts one is trying to get across.  At the same time, one wants to use the ``shiny'' technology to help motivate the students (and, to be honest, the instructor).  By using current software, and using tools that the student can be confident will be used elsewhere, the expense of time needed in explaining the platform can be amortized.  However, some classic tools, such as Maple, can still be very much worth the effort, as we will demonstrate by examples.

This represents both an opportunity and a challenge for the instructor.  There are enough new features in Maple that are useful in instruction that for many of us, some upgrading of our training (or even retraining!) is necessary.  One might question the effort needed, especially given the rise of Python, and the rise of Julia; but we think it is very worthwhile.

What about Mathematica?  We have heard that it, too, has added significant and interesting new features; but we work in an institution with a Maple site license, and in an institution that uses SAGE Math; we will let our comments on Maple stand for (what we believe will be) similar comments for Mathematica.

A bigger question is what about GeoGebra (or Desmos, or Maple Learn).  Our approach has largely been algebraic, or programmatic, and not explicitly visual, although we do use visualization of algebraic and functional concepts.  We do not weigh here the comparative merits of this approach versus the more visual approach of the geometrically-flavoured teaching tools just mentioned; for the sake of argument, we ask the reader to accept (if only temporarily) that an algebraic and programmatic approach can be useful.  We acknowledge that the geometric approach can also be useful, of course: see for example~\cite{Botana2015} or~\cite{tran:2014}.

So: the basic point of view of this paper---which is shared by a great many other papers, including many at this ATCM and at previous ones---is that instructors can, and should, use current technology to help teach mathematics.  In particular, we advocate using \textsl{computational discovery} with the help of the technology to strengthen students' knowledge of \textsl{both} mathematics and of computational technology.  We'll now give specific (algebraic or programmatic) examples of how to do this, and discuss our experiences in the classroom and out.

The literature on using technology for teaching is very, very large.  We will try to put our efforts in perspective in that literature in the final section of this paper, ``Concluding remarks and Further Reading.'' We will include technical references in each specific section.

\section{Specific goals}
\subsection{Course objectives for the students}
\begin{enumerate}[topsep=2pt, itemsep=-4pt]
 \item  To acquire facility in using the computer as a tool for solving and
 exploring mathematical problems.
 \item To learn the fundamental concepts and techniques of procedural and object-oriented
programming. These include: flow control, modular construction, elementary data
structures, recursion, and graphics
\item  To develop problem-solving and communication skills by solving programming projects.
\item To learn to ask new questions, and to be a bit more comfortable with not being able to find out a definitive answer.
\item To become familiar with other computational resources for mathematics,
in particular, but not limited to: \LaTeX, the Online Encyclopedia of Integer Sequences (OEIS), and the Inverse Symbolic Calculator (``identify" in Maple).
\end{enumerate}
\subsection{Learning Outcomes}
On completing this course, the student will be expected to be able to:
\begin{enumerate}[topsep=2pt, itemsep=-4pt]
\item Take a mathematical question and write it as a computational question.
\item Give examples of student-generated questions, including ``what if?'' questions
\item  Give a pseudo-code version of an algorithm to solve the computational question
\item Convert pseudo-code into Sage or Maple or Python code
\item Comment code cleanly and effectively
\item Obtain and interpret output from their code
\item Give visualization and other methods of viewing output effectively.
\item Write up lab reports  in \LaTeX\ or Markdown on their investigations
\end{enumerate}

\subsection{Assessment}
\noindent The course that depends on the resources we describe here is largely project based: there typically are small projects
(especially at the beginning of the course) to ensure that everyone is coming
up to speed with coding.  Then there can be several more substantial projects,
taking a mathematical question, and learning how to explore it in an
experimental fashion.  We make no assumptions about level of computer
experience: if the student has written code before that will help, but if they
haven't we will get them up to speed quickly.

Their projects are assessed using the following criteria:
\begin{itemize}[topsep=2pt, itemsep=-4pt]
 \item Does the project contain new\footnote{By this we mean \emph{new to the student}, ideally questions that they came up with themselves.  If they are genuinely open questions, or genuinely novel questions not in the literature, so much the better; but this is not required.} questions, new thoughts, or new answers?
 \item Do they have code that {\sl runs}?
 \item Does the code produce {\sl output}?
 \item Is the output {\sl correct}?
 \item Is the output {\sl complete}?
 \item Is the code {\sl clean}?  Are the variables suitably named?  Are
 functions, control structures, data types, etc used appropriately?
 \item Is the code {\sl documented}?
 \item Is the project written up as a lab report?
\end{itemize}
If the project passes all these criteria, it will get a good grade!

In some situations (e.g. for oral presentations) \emph{peer-assessment} may be used: that is, part of your grade will be determined by fellow students.  Each student will be expected to supply assessments of their fellow students in turn\footnote{Telling the students that the peer-assessed marks will be blended with the instructors' assessment, because students are sometimes harsh on each other, tends to improve things.  We also ask students to justify their assessments.}.
\section{The first example}
Consider the sequence
\begin{align}
    1, \frac32, \frac{17}{12}, \frac{577}{408}, {\frac{665857}{470832}}, \ldots\>.
\end{align}
Plunked down in front of the students without explanation, it is pretty mysterious.  Explicitly telling them to square each entry generates a more intelligible sequence, namely
\begin{align}
    1, \>\frac94 = 2\frac{1}{4} = 2\frac{1}{2^2},\> \frac{289}{144} = 2\frac{1}{12^2},\> {\frac{332929}{166464}} = 2\frac{1}{408^2},\> {\frac{443365544449}{221682772224}} = 2\frac1{470832^2}, \ldots\>.
\end{align}
Although more intelligible, this sequence is \emph{still} mysterious.  At this early point in the course, students are reluctant to ask questions (not just for fear of looking foolish in the eyes of their peers, but also for fear of being wrong~\cite{holt1982children}), but the urge to ask is, for many, irresistible here.  We get questions like, ``where did this sequence come from?'' and ``is the sequence of squares really going to $2$?'' and, if we are lucky, ``does this mean that the original sequence tends to $\sqrt{2}$?''
When we ask the Online Encyclopedia of Integer Sequences (OEIS) about the sequence of (say) denominators of the original sequence, we get~\cite[A051009]{sloane2021}, the ``Reduced denominators of Newton's method for sqrt(2),'' and at that point we can explain that the rule we were using to generate the sequence was
\begin{equation}
    x_{n+1} = \frac12\left( x_n + \frac2{x_n}\right)\>;
\end{equation}
that is, we start with an initial estimate ($x_0=1$), and then we average it with what we get when we divide it into $2$.  If our initial estimate of the square root had been correct, then we would stay the same; if our estimate was too small, then dividing it into $2$ would give a larger number and averaging is plausibly a way to generate a better answer.  It's still mysterious to the students, but they are beginning to get interested.

It gets better: we then expand into continued fractions.  For space reasons we use a compact notation here.  The original sequence is $1$, $1+[2]$, $1 + [2,2,2]$, $1 + [2,2,2,2,2,2,2]$, and so on; the students are astonished, perhaps a little appalled.  We do the first few by hand, and then the rest by using Maple's ``convert'' facility (to ``confrac'').

We can then digress to the other continued fractions $1+[2,2,\ldots,2]$ where the number of $2$'s is not one of $1$, $3$, $7$, $15$, and so on.  This is already a very fruitful line of investigation.

Indeed most students are stunned to discover that the square root of $2$ has a perfectly predictable infinite continued fraction.  Getting them to make the claim that
\begin{equation}
    \sqrt{2} = 1 + [2, 2, 2, \ldots ]\>,
\end{equation}
even given the numerical evidence above, takes encouragement from the instructor. The students are afraid, at this point, of a trap, of guessing wrong, especially when they ``know'' (i.e., have been taught) that the decimal expansion of the square root of two has no pattern.  This naturally leads to a discussion what it means to perform an infinite number of operations.  The OEIS proves useful here, as well; perhaps simply by letting the students in on the infinitely wonderful world of sequences.

From here, we can go to true Newton's method, in general, and not just for square roots.

\section{Visualizing Newton's Method}
\begin{figure}
    \centering
     \begin{subfigure}[b]{0.45\textwidth}
         \centering
         \includegraphics[width=\textwidth]{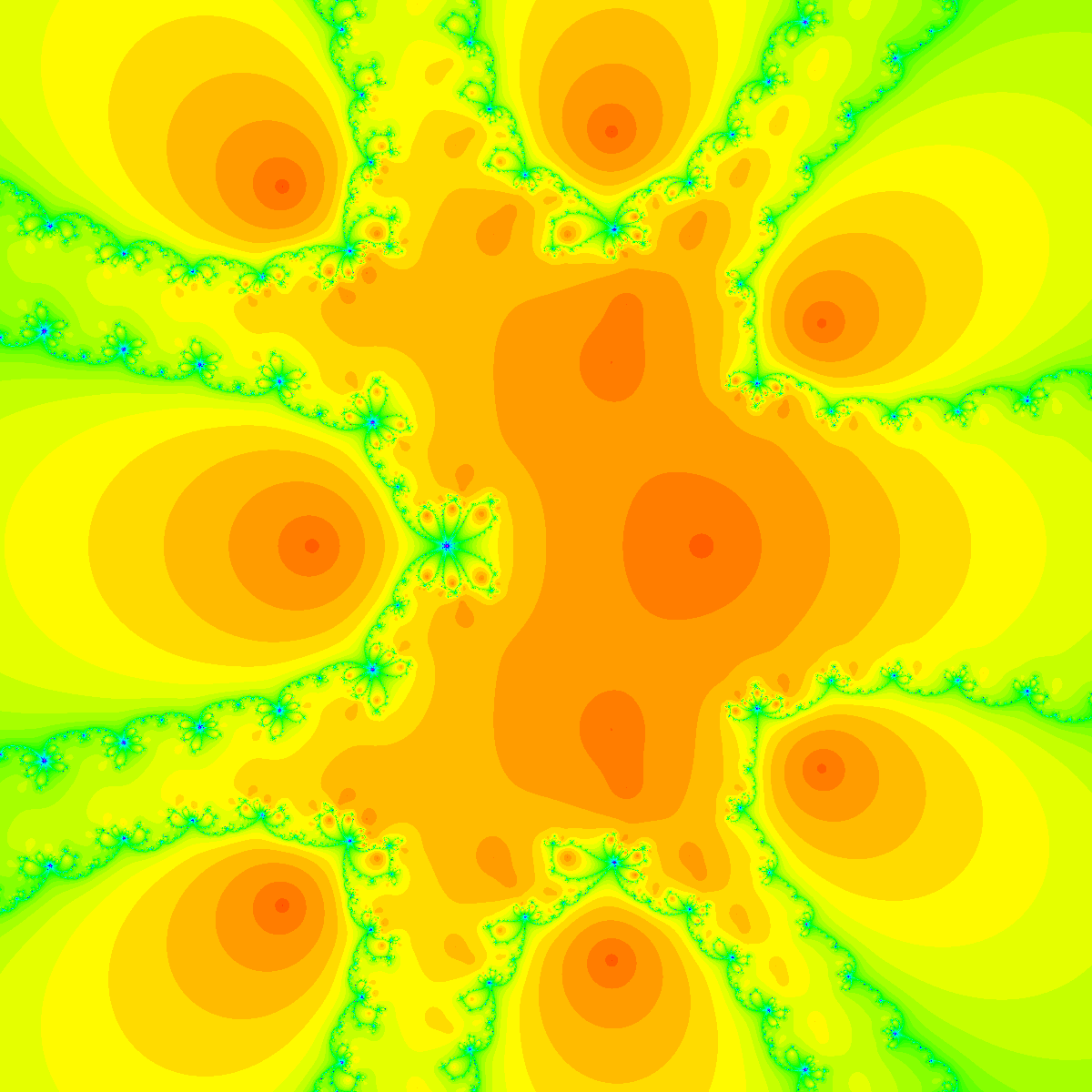}
         \caption{Maple version}
         \label{fig:fibmand6Maple}
     \end{subfigure}
     \hfill
     \begin{subfigure}[b]{0.45\textwidth}
         \centering
     \includegraphics[width=\textwidth]{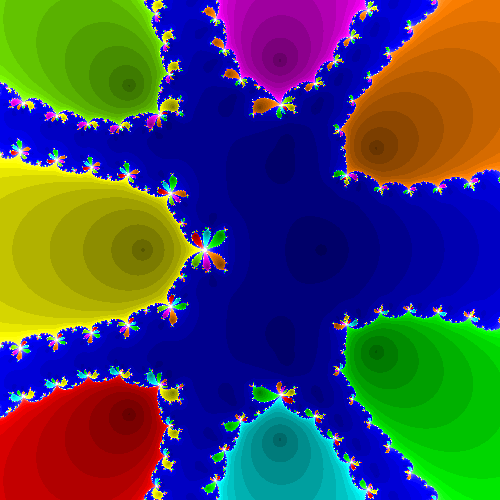}
     \caption{Python version}
     \label{fig:fibmand6Python}
     \end{subfigure}
    \caption{Basins of attraction for Newton's method on a certain degree~$8$ polynomial, displayed in equation~\eqref{eq:FibMand6}, as produced by the \texttt{Fractals:-Newton} command in Maple (left) and in by a custom Python program (right). The view in the complex plane is from $-1.8 \le x \le 1.0$ and $-1.4 \le y \le 1.4$.}
    \label{fig:FibMand9}
\end{figure}

The chaotic dynamics of Newton's method are very well-studied: see for instance~\cite{Vrscay1986}, \cite{Strang1991}, \cite{Corless1998}, \cite{li2019revisiting}, \cite{Calkin:2021:Newton}, or indeed~\cite{mandelbrot2002some}. They prove perennially popular with students, in part because they produce such pretty pictures.  See Figure~\ref{fig:FibMand9}, which shows the basins of attraction for Newton's method on the polynomial
\begin{align}
    z^{8}+3 z^{7}+5 z^{6}+5 z^{5}+4 z^{4}+2 z^{3}+z^{2}+z\>. \label{eq:FibMand6}
\end{align}

We want the students to learn to produce their own such pictures, and to begin to understand what they mean.  To this end, we use both built-in procedures and hand-crafted procedures. The best ones of all are the ones the students learn to write for themselves, but we have to work up to this.

\subsection{Summary of what we teach in this section}
We teach the students
\begin{enumerate}[topsep=2pt, itemsep=-4pt]
    \item The special Newton iteration $z_{n+1} = (z_n + a/z_n)/2$ for finding square roots.
    \item The general Newton iteration formula $z_{n+1} = z_n - f(z_n)/f'(z_n)$; we usually give a graphical demonstration similar to the usual one sketched in a first calculus class.
    \item What happens when you apply Newton iteration to $f(z) = z^2 - 2 $, using exact rational arithmetic and $z_0 = 1$, and its connection to continued fractions (as previously discussed).
    \item A game, which we call ``pass the parcel,'' which cements the idea of iteration: the students must take turns evaluating the Newton iteration and then pass their answer on to the next group to compute the next iterate.
    \item The importance of the initial estimate $z_0$ (more on this, below, but we try to get the students to realize this without us telling them; this actually happens, especially after playing ``pass the parcel,'' and the students do ask how to choose the initial estimate).
    \item That $z_n$ exactly solves the \textsl{nearby} equation $f(z) - f(z_n) = 0$ and show that $f(z_n)$ gets very small very rapidly if the initial estimate is good.  We interpret this as a change in the question: e.g. $17/12$ is not just an approximation to the square root of two, it is the exact square root of $2\,\sfrac{1}{144}$; that is, it is the exact square root of an approximation of two, and we can compute this \textsl{residual error} even if we don't know the true value of the square root of two.
    \item The very beginnings of floating-point arithmetic and its peculiarities.  This is actually very hard to teach at this level, but we believe it is necessary.  We also teach the beginnings of complexity theory, by showing that exact arithmetic is costly; floating-point is vastly faster, \emph{because its memory usage is predictable}.
    \item How to use the Online Encyclopedia of Integer Sequences (OEIS) at \href{http://oeis.org/}{http://oeis.org/}.
    \item What happens when you apply Newton iteration with a \textsl{real} initial guess to $f(z) = z^2+1$; see~\cite{Strang1991} and~\cite{li2019revisiting}.
    \item How to use the Maple  \texttt{Fractals} package to make fractals of their own.
    \item Variations on Newton iteration including the secant method, Halley's method, Householder iteration, and Schr\"oder iteration; and that all of these (apart from the secant method) can be considered as Newton iterations of \emph{some} function (although the proof of that requires calculus).
\end{enumerate}
\subsection{Behind the curtain: why we teach what we do, and what the students get out of it}
The main pedagogical purpose of this section is to reinforce the notion of a \textsl{function}.  Students usually have been taught functions and even differentiation rules in high school, but frequently they differentiate only by rote and their notion of a function is weak: they typically only think of functions as expressions.  By getting them to play ``pass the parcel'' (with the support of a Maple operator or procedure for the function) we emphasize the \textsl{active} nature of a function.  This also is an easy way to break the ``inactivity barrier'' because they have to interact with their fellow students.

We use \emph{scaffolding}: we supply the functions $f(z)=0$ to solve and the initial estimates $z_0$, at first; but then we give some functions \emph{without} initial estimates.  The students quite rightly find it difficult to construct their own (although some ingenuity is often displayed) and are naturally led to ask the crucial question of what is the influence of the initial estimate.  This leads directly to the Newton fractal pictures.

The secondary pedagogical purpose is a gentle introduction to the software tools.  We also use Python for this section, but in that case there is more programming involved\footnote{What the students gain from this increase in programming effort includes greater control: for instance, the picture in Figure~\ref{fig:fibmand6Python} includes more intelligible detail, in that different colours correspond to different limits in the iteration, whereas the colours in the Maple fractal were chosen generically by the developer (David W.~Linder at Maplesoft) to indicate iteration count only.}, up to and including \textsl{automatic differentiation} (as opposed to \textsl{symbolic} differentiation, or as opposed to differentiation by hand, or as opposed to numerical differentiation by finite differences or other numerical techniques).

A third pedagogical point is to introduce the beginnings of complexity theory---computation with exact rationals is expensive and the growth of the length of the exact rational answers is remarkable; students are quickly converted to the worth and utility of floating-point arithmetic.  Newton's method is relatively benign for use in floating-point, usually; this is a help.  We only introduce as much numerical analysis of floating-point numbers as we are forced to do.

There are many directions to take this unit: we make connections to number theory via the theory of continued fractions, for instance.  We choose to do this in part because continued fractions are not part of (most) students' curricular choices, and therefore we are not ``stealing the thunder'' from subsequent courses\footnote{This curious English expression means, in this context, to use the highlights; ``stealing thunder" would make this particular course more exciting, but would make subsequent courses more boring, which is unfair to the downstream teacher.}. The main reason, though, is that the results are surprising and beautiful.

Using the OEIS, for instance, the students can be guided to discover that the $n$th iterate of Newton's method for $f(z) = z^2-m$, starting at the rational initial estimate $z_0 = a/b$, is
\begin{equation}
    z_n = \sqrt{m}\left(\frac{1 + r^{2^n}}{1 - r^{2^n}} \right)
\end{equation}
where $r = (a-b\sqrt{m})/(a+b\sqrt{m})$. Since $0 < r < 1$ if $a$, $b$ and $m$ are positive, this shows the quadratic convergence of Newton's method very well\footnote{This purely algebraic formula is equivalent to the cotangent formula of~\cite{Strang1991} and~\cite{li2019revisiting}.}.
There are other number-theoretic directions to go, as well: one can use this as a springboard to study Pell's equation.  There are higher-dimensional versions, and indeed the \emph{matrix} square root is a subject of much research---see for example~\cite{higham1986newton} and \cite{Higham1997}.

We try to avoid using the word ``convergence,'' however, and satisfy ourselves by computing the exact square roots of numbers near to $2$ (for instance).  This is a surprisingly useful intermediate step.

\subsection{Suggestions for Assessment}
It is straightforward to ask the students to compute Newton iterations for other quadratic irrationals, and for other algebraic numbers.  It is also straightforward to get them to use Maple's built-in Newton fractal package to find basins of attraction for functions of their own choice.  Our most successful exercise here so far, however, has been to hold a contest to see who could find a function that produced the ``best'' picture (by choosing their own functions and regions). Unsurprisingly, the students overwhelmingly produce excellent pictures, and we have sometimes settled on making a collage of them all.  The students are then asked to make a similar collage of the basins of attraction of Halley's method\footnote{In some classes we covered quite a few iterative methods: see~\cite{li2019revisiting} for a case where a student solved a problem (in class!) that was believed to be open.} on the \emph{same} functions; by asking Maple to do Newton iteration on $f(x)/\sqrt{f'(x)}$ one is actually asking for Halley iteration\footnote{This fact is in the literature, as is the more general fact previously mentioned that \emph{every} iteration can be cast as a Newton iteration for \emph{some} function.}.  The resulting pair of collages is quite instructive in that one can understand some differences between Newton iteration and Halley iteration intuitively thereby.

\section{Bohemian Matrices}
\begin{figure}
    \centering
    \includegraphics[width=0.5\textwidth]{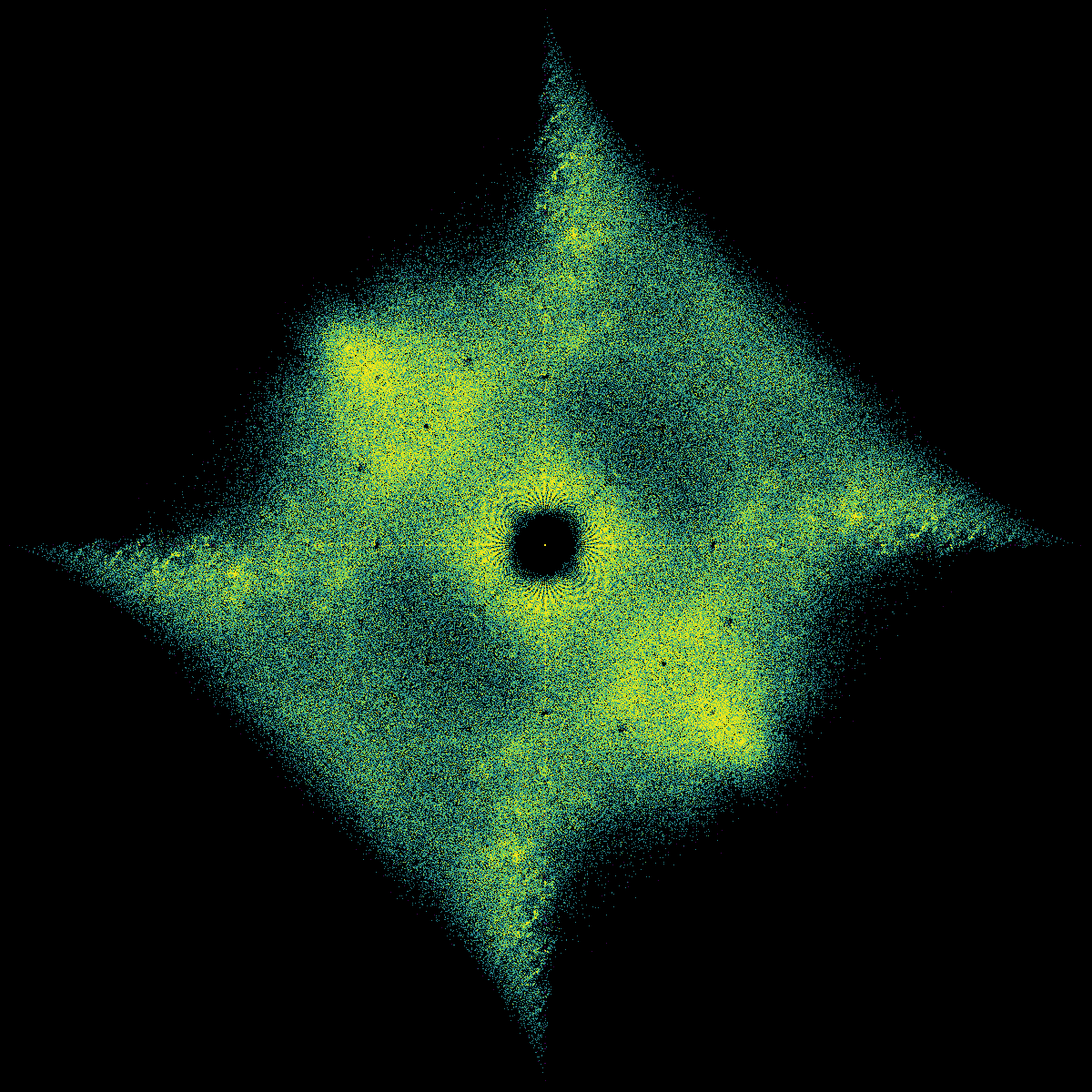}
    \caption{Eigenvalue density plot of all $2^{17} = 131072$ ten by ten Bohemian skew-pentadiagonal matrices with population $\{1,i\}$, plotted on $-3.25 \le \Re(\lambda) \le 3.25$ and $-3.25 \le \Im(\lambda) \le 3.25$ in the complex $\lambda$-plane.  Hotter colours have higher density.  At this time of writing we have no explanation whatever for the distribution pattern visible here.}
    \label{fig:NewBohemian}
\end{figure}

The section of the course on \textbf{Bohemian Matrices} is popular, in part because there are so many interesting pictures already, but mostly because there is so much unknown.  Many of the pictures in the gallery at \href{www.bohemianmatrices.com}{www.bohemianmatrices.com} were actually produced by students in the course, and never seen by anyone before those students made them.  In Figure~\ref{fig:NewBohemian} you can see a new picture made specifically for this paper\footnote{This uses \emph{skew-symmetric} \emph{pentadiagonal} matrices; the five by five case with generic entries is
\begin{equation*}
    \left[\begin{array}{ccccc}
0 & t_{1} & t_{5} & 0 & 0
\\
 -t_{1} & 0 & t_{2} & t_{6} & 0
\\
 -t_{5} & -t_{2} & 0 & t_{3} & t_{7}
\\
 0 & -t_{6} & -t_{3} & 0 & t_{4}
\\
 0 & 0 & -t_{7} & -t_{4} & 0
\end{array}\right]\>.
\end{equation*}
Each $t_i$ must be chosen from the finite population; all possible choices make up the Bohemian family.  Here if $P= \{1,i\}$ there are $2^7 = 128$ possible matrices in this family with this dimension.
}, and likewise never seen before by anyone (there are infinitely many such; welcome to the party).

\subsection{Summary of what we teach}
We don't wish to steal the thunder of linear algebra classes, so we teach eigenvalues from a determinantal point of view.  We begin with the single linear equation $t_1 x_1 = 1$ in the single unknown $x_1$ where $t_1$ is known to be drawn from our finite population (say, $t_1$ must be a member of $P = \{-1, 0, 1\}$).  Then, unless $t_1=0$, we can solve this system of equations; that is, whether $t_1=0$ or not \emph{determines} whether or not this equation can be solved for $x_1$.
We then move on to a related two-by-two system of equations in the two unknowns $x_1$ and $x_2$:
\begin{align}
    t_1 x_1 + t_2 x_2 &= 1 \nonumber\\
    - x_1 + t_1 x_2 &= 0
\end{align}
Only one new member of the population, $t_2$, is involved.  By multiplying the second equation by $t_1$ and adding the first equation, we get $(t_1^2 +t_2)x_2 = 1$, and now $t_1^2 + t_2$ determines whether or not we can solve the equations for $x_1$ and $x_2$.

The three-by-three case is analogous:
\begin{align}
    t_1 x_1 + t_2 x_2 + t_3 x_3 &= 1 \nonumber\\
    - x_1 + t_1 x_2 + t_2 x_3 &= 0 \nonumber\\
            -x_2 + t_1 x_3 &= 0
\end{align}
and the students are easily guided to discovering that the determinant of this (Toeplitz) system is $t_1^3 + 2t_1t_2 + t_3$.  The determinant of the $4$ by $4$ version is
$t_{1}^{4}+3 t_{1}^{2} t_{2}+2 t_{1} t_{3}+t_{2}^{2}+t_{4}$, and the $5$ by $5$ version has determinant $t_{1}^{5}+4 t_{1}^{3} t_{2}+3 t_{1}^{2} t_{3}+3 t_{1} t_{2}^{2}+2 t_{1} t_{4}+2 t_{2} t_{3}+t_{5}$.  The students have an interesting time trying to guess the pattern here, and even experts can be stumped: but we will leave it as an exercise.

Replacing $t_1$ with $t_1 - \lambda$ is straightforward; one then has a collection of polynomials to solve, which the students are happy to leave to the machine.  One can then show them an ``eigenvalue'' solver without shocking them too much, or doing too much violence to the later curriculum.
\subsection{Why we teach it}
Eigenvalues are often given short shrift in a first linear algebra course, and the students need more practice with them.  We can also use the occasion to show the students different matrix structures including symmetric matrices, skew-symmetric matrices, and others. The structure of the introductory example is a special case of a Toeplitz matrix structure.   We even get the students to invent their own: the ``checkerboard'' matrix picture at www.bohemianmatrices.com is an example of a student-generated matrix structure, which we've never seen elsewhere and seems unlikely to have any application, but we like it anyway.

But the biggest reason to have this as part of the course is that the research area is so new that even many easy questions are as yet unanswered, and the students have quite a good chance to contribute something new.
Frankly, since there are many open problems here, and we are just scratching the surface of this new field, we can use the students' creativity, too.
\subsection{Suggestions for Assessment}
Asking students to count quantities of interest and to try to find patterns can be quite accessible, for some matrix structures.  For instance, one can ask how many $m$ by $m$ Toeplitz matrices of the form used above there are, with population $P = \{-1,0,1\}$ (this is easy); one can ask how many distinct characteristic polynomials there are at dimension $m$; one can ask how many different eigenvalues there are at dimension $m$ (this is quite a bit harder, because some eigenvalues may be shared between more than one matrix in the Bohemian family); one can ask how many multiple eigenvalues there are, or how many matrices are singular, and so on.  One can ask \emph{programming} questions: how quickly can you generate all $3^m$ matrices, for instance.

Our favourite assessment here, however, came in the form of a contest: who could produce the most interesting pictures?  As a ``filtering'' mechanism, this failed utterly: the students all loved each others' pictures, and refused, quite rightly, to rank them.  We were quite happy to give everyone full marks for this exercise.
\section{The Mandelbrot Polynomials \label{sec:Mandelbrot}}
For general material on the \emph{Mandelbrot set}, the
\href{https://en.wikipedia.org/wiki/Mandelbrot_set}{Wikipedia page} on the subject provides a good starting point to learn more.
Indeed many students are familiar with the definition of the {Mandelbrot set}, but few of them have seen the \emph{Mandelbrot polynomials}, defined by
$z_{n+1} = z_n^2 + c$ where $z_0 = 0$ and $c$ is a symbol or variable representing an as yet unknown complex number.  The first few of these polynomials are listed in~\cite{Calkin:2021:Mandelbrot}: $0$, $c$, $c^2 + c$, $c^{4}+2 c^{3}+c^{2}+c$, and so on.

By studying the \emph{zeros} of these polynomials we are studying the \emph{periodic points} of the Mandelbrot iteration.  This turns out to be surprisingly fruitful, in that there is a connection to eigenvalue problems, as detailed in the aforementioned paper and its references, as well as a connection to the modern theory of dynamical systems.
\subsection{Summary of what we teach}
We use this example to teach that computing roots of polynomials can be numerically difficult: the condition number $B_n(c) = \sum_{k=0}^{d_n} |a_k| |c|^k$ of the Mandelbrot polynomials $z_n(c) = \sum_{k=0}^{d_n} a_k c^k$ grows \emph{doubly exponentially} with $n$ (singly exponentially with the degree $d_n = 2^{n-1}$) which means that the number of decimal digits needed to compute the zeros grows linearly with $n$.  In contrast, the matrix eigenvalue problem is well-conditioned, and can be used in ordinary double precision floating-point arithmetic for remarkably large dimension of the matrix and degree of the polynomial.
\begin{figure}
    \centering
    \includegraphics[width=8cm]{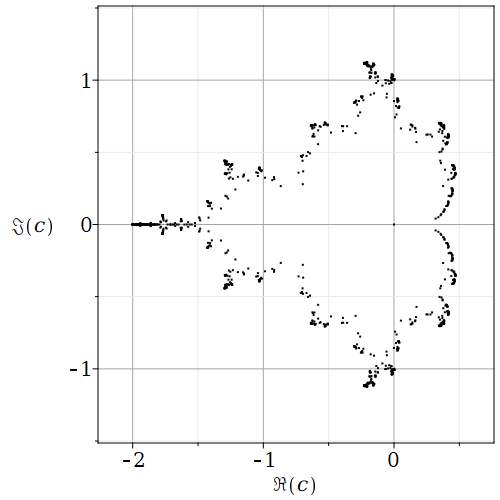}
    \caption{Zeros of $z_{11}(c)$, a polynomial of degree $1024$. Each of these zeros is a value of $c$ which leads to a periodic orbit of period $11$, namely $0$, $z_1(c) = c$, $z_2(c)$, $\ldots$, $z_{10}(c)$, $0$.  Such points are also called \emph{hyperbolic centers} in the Mandelbrot set and are actually quite sparse in the Mandelbrot set.}
    \label{fig:MandelbrotZ11}
\end{figure}

\subsection{Why we teach it}
The Mandelbrot polynomials give a good excuse for studying \emph{nonlinear} iterations in general; it is a very great surprise that, for $c$ not in the Mandelbrot set, this recurrence relation can be ``solved'' exactly (well, by a convergent series); see~\cite{Calkin:2021:Mandelbrot} and~\cite{MichaelLarsen2020}.  We did not know this ourselves when we were teaching from this material, and we now wonder just how much of this we can include in a first course.  But Peter Taylor of Queen's University\footnote{This is from a talk he gave at the Summer Canadian Math Society Meeting 2021.} has been able to successfully use a similar iteration and its asymptotics to study a bacterial growth model in a first-year context, and so we hope to be able to do so as well.

\subsection{Suggestions for Assessment}
Students enjoy writing fast programs to draw the Mandelbrot set; getting them to carefully draw pictures of the periodic points in the Mandelbrot set is not difficult.

There are $21$ ``facts'' listed in~\cite[Sec. 2]{Calkin:2021:Mandelbrot}, some of which make good questions for students to attempt; in the following section there are seven \emph{conjectures}, open as we write this.  The last two, asking if Mandelbrot polynomials are \emph{unimodal}\footnote{A polynomial is unimodal if its nonzero monomial basis coefficients are all positive integers, and those coefficients increase towards a unique maximum and then decrease~\cite{wilf2005generatingfunctionology}.}, seem as if they ought to be accessible.  However, we don't know how to prove either, and one of us tried quite hard.  Perhaps a student could do it. If so, they would get a paper out of it.

Getting them to \emph{come up with their own questions}, however, is one of the goals of the course.  By this time in the course, we hope that the students are doing so: perhaps wondering what happens if the constant $c$ is replaced by a pre-specified sequence, or a random sequence; or what happens if $c$ is not a scalar but a matrix; or if the iteration is $z_{n+1} = z_n^\alpha + c$ for some noninteger $\alpha$, possibly even complex~$\alpha$.  Some of these questions are sterile, as will be some of the questions the students come up with; but some might not be, and the adventure begins.
\section{Concluding Remarks and Further Reading}
Our aim is to try to teach an experimental approach to mathematics, \emph{without} either stealing the thunder of courses already in the standard curriculum, or requiring significant prerequisites.
Mathematics is vast and the standard undergraduate curriculum only scratches its surface.  Nonetheless, those scratches are somehow \emph{fractal} and it is difficult to avoid them while remaining relevant and interesting.   As presented here, this material comes right up to and touches on several fundamental notions of the standard mathematics curriculum: for example, we want the students to \emph{experience} convergence, which motivates the theory of limits. We also touch on existence and uniqueness of solutions of equations, and the meaning of proof.  We believe that this extra motivation for standard material is to the good.

The use of advanced mathematics
to teach early undergraduate students was advocated by Mandelbrot himself. See for instance~\cite{mandelbrot2013fractals}, \cite{frame2002fractals}, and~\cite{mandelbrot2002some}.

We are not trying to get the students to a destination: we want them to experience the journey.  It is less like using a GPS, and
more like an extended Sunday exploration,
perhaps with a paper map.

When introducing students to the ideas of experimental mathematics, learning how to ask questions is not just
a part of the process, it {\sl is } the process.

Therefore, we concentrate much more heavily on things \emph{not} in the standard curriculum: such as actual mathematically-oriented computer programming, which strengthens the notion of a function. This kind of programming also both motivates, and partially replaces, the notion of proof, for this course. A working program helps you check that you have an effective definition.

Other authors have tried active non-standard approaches before, of course.  One that comes to mind is Z.~A.~Melzak's \emph{Companion to Concrete Mathematics}~\cite{melzak2007companion}, which is a rich source of nonstandard tricks and the problems they apply to. The most venerable of such works, however, is the magnificent two-volume set by P{\'o}lya and Szeg\H{o}~\cite{polya1997problemsII,polya1997problemsI}.

Why use active learning at all?  The evidence for its effectiveness now is so strong---see, for example, \cite{freeman2014active}---that it is unethical \emph{not} to use it, if you can.  It is as simple (and as difficult) as that.  The mathematics education literature is very nearly as vast as the mathematics literature\footnote{In addition to the many volumes of ATCM proceedings, see, for instance, the chapters in~\cite{ME4AI:2021} and the references therein, as a slice of the intersection of the math ed literature with AI literature.}, but there are few ideas so well-supported in that literature as the now-established fact that it is better for students to \emph{do mathematics} than it is for them to listen or watch it and then regurgitate it on an exam. Indeed, quoting from~\cite{yang2013discovering}, ``Creativity does not come from drills.''

There are principles that help the teacher design good exercises to increase the activity of the student: a good exercise should be connected to great currents of mathematical thought (so as not to waste the student's time); it should be engaging, perhaps by having an element of surprise about it; it should be accessible at the student's level of education; it should be rich and allow open-ended exploration.  Those two last criteria are sometimes termed ``low floor, high ceiling,'' meaning the exercise should be easy to start, but have many further levels.  We hope that you agree that the examples presented here fit these criteria.

Where the literature diverges, however, is on how best to actually \emph{get} the students to be active.  For example, the cumulative online proceedings of this conference since at least 1997 show a very wide variety of approaches, and it is quite possible that they \emph{all} work, perhaps to varying degrees.
What we are recommending here, like many others before us, such as~\cite{abelson1986turtle}, is to \emph{enrich the content}~\cite{corless2004computer}; as in~\cite{betteridge2019teaching} and~\cite{betteridge2020teaching} we are describing courses specifically designed to encourage active learning of mathematics.  Given that many curricula are already overfull, the question of where this will fit in the student's timetable is a difficult one, but it seems clear that to make room for this course, some other less-useful course will need to be pushed aside.  We, innocently, will look the other way, whistling, as you decide just which.

\textbf{Acknowledgements.} We thank Professor Wei-Chi Yang for the invitation to present this work, and the anonymous referees for very helpful remarks.  This work was partially supported by an award from the Centre for Teaching and Learning at Western University.

\bibliographystyle{plain}


\end{document}